\documentclass[10pt,letterpaper,twocolumn]{article} 

\usepackage{ol2}
\usepackage[draft]{hyperref}
\usepackage{amsmath}

\usepackage{epstopdf}

\begin{document}

\twocolumn[ 

\title{An elementary approach for the phase retrieval problem}


\author{Yuan Sun$^{*}$}

\address{
Department of Physics and Astronomy, State University of New York at Stony Brook\\ Stony Brook, NY 11794-3800, USA
\\
$^*$Corresponding author: yuan.sun.1@stonybrook.edu
}

\begin{abstract}If the phase retrieval problem can be solved by a method similar to that of solving a system of linear equations under the context of FFT, the time complexity of computer based phase retrieval algorithm would be reduced. Here I present such a method which is recursive but highly non-linear in nature, based on a close look at the Fourier spectrum of the square of the function norm. In a one dimensional problem it takes $O(N^2)$ steps of calculation to recover the phases of an N component complex vector. This method could work in 1, 2 or even higher dimensional finite Fourier analysis without changes in the behavior of time complexity. For one dimensional problem the performance of an algorithm based on this method is shown, where the limitations are discussed too, especially when subject to random noises which contains significant high frequency components.
\end{abstract}


 ] 

\noindent A fast algorithm for the phase retrieval problem could be vital to many fields of optics and related areas, such as X-ray imaging, electron scattering, and the phase reconstruction of a laser pulse. Up to today, the main stream algorithms are those based on the Gerchberg-Saxton algorithm\cite{G-S}, especially the hybrid input-out method \cite{Fienup82}, which are all iterative in nature. There are also many pioneering works on non-iterative methods, such as those based on the fractional Fourier transform\cite{Cong98} and those based on `locating the zeros'\cite{Deighton}. A lot of deep discussions are available in the literature too about the uniqueness and existence of the solution, of interest is a recent paper\cite{Candes}, yet those works are a little abstract to some extent.

In the last decade, many interesting papers have been published on this subject. \cite{PhysRevA.75.043805} experimentally showed a simple and robust approach for retrieving arbitrary complex-valued fields from three or more diffraction intensity recordings. And \cite{Marchesini07} provides a thorough evaluation and comparison of iterative projection algorithms for phase retrieval.

The phase retrieval problem is highly nonlinear in nature. Hence attempts to convert it into a linear problem may not be necessarily effective. In this sense iterative methods do seem to be appropriate, except for the minor issues of convergence, especially the problem of stagnation\cite{Fienup86}. However, we would still wonder, whether an elementary approach exists, just like that of solving a system of linear equations?

Dallas' algorithm\cite{Dallas75} opens an interesting new direction to the answers to the above problems. The original Dallas' ideas, together with many following works \cite{Fienup83} \cite{Crimmins87}, have laid the ground of these types of algorithm. Yet, two shortcoming generally exists except for a few special examples: they are well known to be highly sensitive to noise and generate large numbers of alternative branches which are very hard to trace.

Here I wish to present an algorithm which could be viewed as a variation of Dallas' algorithm, which, at the same time, is equipped with an add-on to avoid the two glaring shortcomings.

Let us begin with a very simple observation. Suppose we have $f = \sin \theta$, then ${|f|}^2 = \sin^2 \theta = \frac{1 - \cos 2\theta}{2}$. It seems the square of the norm of the function would incur a Fourier component of higher frequency. This is no coincidence. To see this, I first want to formulate the easiest (but not necessarily trivial) case of the phase retrieval problem under the context of Fourier series, namely, the trigonometric polynomials: $f(x) = \sum_{n = - N}^{N} a_{n} e^{inx}$. Suppose that $| f |$ and all $|a_n|$'s are known. The question is to recover $f$'s or all the $a_n$'s phases. Notice that $|f|^2= \sum_{n = - N}^{N} a_{n} e^{inx} \times \sum_{n = - N}^{N} a^{\dagger}_{n} e^{-inx} \nonumber $ and let us rearrange its terms as in \eqref{f^2series}.

\begin{equation}
\label{f^2series}
|f|^2 = \sum_{l=-2N}^{2N} (\sum_{j= max(l-N,-N)}^{min(l+N,N)}a_j a^{\dagger}_{j-l}) e^{ilx}
\end{equation}

If we expand $|f|^2$ into Fourier series $|f|^2 = \sum_{l = - 2N}^{2N}  b_{l} e^{ilx}$ then by \eqref{f^2series} we see that $b_l = \sum_{j= max(l-N,-N)}^{min(l+N,N)}a_j a^{\dagger}_{j-l}$.  In the context of Fourier transform this is nothing more than saying $\it{FT}(ff^{\dagger}) = \it{FT}(f)*\it{FT}(f^{\dagger})$. The highest order term is that of $e^{i2Nx}$, and the coefficient is $b_{2N} = a_N \cdot a_{-N}^\dagger$. All the $b_l$'s can be calculated from the given $|f|$. Because we know the norm of $a_N$ and $a_{-N}^\dagger$, if we assign a phase to $a_N$, which can be regarded as nothing more than an overall phase we always have the freedom to choose, then we acquire the phase of $a_{-N}$. Next let us move on to look at the coefficient of the next order: $b_{2N-1} = a_N \cdot a_{-(N-1)}^{\dagger} + a_{N-1} \cdot a_{-N}^{\dagger}$. Here, in order to solve for the phases of $a_{-(N-1)}$ and $a_{N-1}$, we are forced to look at an equation of the following type as in \eqref{vector triangle}, where $x$, $y$ and $z$ are known complex number and we are solving for for $\alpha_1$ and $\alpha_2$ which are angles ranging from $0$ to $2\pi$.

\begin{equation}
\label{vector triangle}
x \cdot e^{i\alpha_1} + y \cdot e^{i\alpha_2} = z
\end{equation}

If the original phase retrieval problem is not ill-defined then of course the above equation \eqref{vector triangle} for $a_{-(N-1)}$ and $a_{N-1}$ shall be satisfied and the existence of a solution is out of the problem. If for some reason \eqref{vector triangle} can not be satisfied we can still choose $\alpha_1$ and $\alpha_2$ such that the left hand side and the right hand side can be as close as possible. The idea of \eqref{vector triangle}, roughly speaking, is to solve for angles in a triangle where the lengths of all three sides are given: $a = |x|$, $b = |y|$ and $c = |z|$. This can be easily done by the law of cosines, but with an important subtlety here: we can have two sets of solutions, as shown in (1) and (2) of Fig 1.

Let us pretend for a moment that we can somehow choose the one right solution out of two possibilities in \eqref{vector triangle} henceforth get $a_{-(N-1)}$ and $a_{N-1}$. Then we can go on to look at $ b_{N-2}$ and will again get a \eqref{vector triangle} type of equation for $a_{-(N-2)}$ and $a_{N-2}$. Say again we can luckily choose the right solution out of two possibilities. We can keep on going to the next order of $b_l$ recursively and solve a \eqref{vector triangle} type of equation repeatedly ($N$ or $N+1$ times) to fully recover the phases of the Fourier coefficients' phases of $f$ -- requires $O(N)$ steps of calculation to recover $2N$ phases -- provided we can always choose the correct one solution out of two in solving \eqref{vector triangle}. This process is recursive, but still very non-linear in nature (as we are solving a lot of triangles instead of a linear system), as compared to the linear nature recursive method of fractional Fourier transform as in \cite{Cong98}.

\begin{figure}[htb]
\label{TriangleGraph}
\centerline{\includegraphics[width=7.5cm]{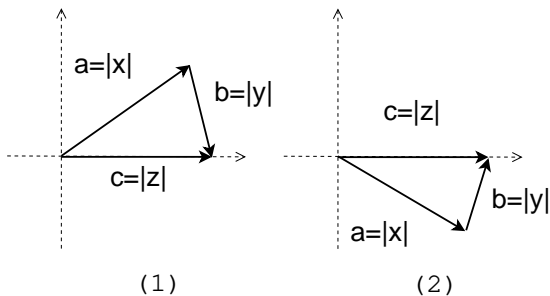}}
\caption{Solving the triangle for equation \eqref{vector triangle}}
\end{figure}

The computer algorithms for Fourier transforms are generally FFT, which is based on the context of finite Fourier analysis. Now let's reformulate a bit to fit this context, to establish the finite Fourier transform between $(0,\cdots,0,a_{-N/2},\cdots,a_{N/2-1},0,\cdots,0)$ and $(F(-(N-1)),...,F(N-1))$, as shown in \eqref{FFAforward} and \eqref{FFAbackward}. Here it is implicitly assumed that $F$ itself is a $2N-1$ component vector while its Fourier transform has only $N$ nonzero components in the lower frequency. This is to make space for the higher frequency part of $|F|^2$ to show up, in analogy to the Fourier series context discussed previously. Theoretically speaking, if both $a_l$ and $F(k)$ are to have N components, then the higher frequencies in $|F|^2$ would not be properly demonstrated by $N$-point finite Fourier analysis; in practice, in order to reveal the higher frequency components in $|F|^2$, we have to double-oversample $F$. If the original phase retrieval problem is not ill-posed this oversampling of $F$ shall carry out naturally.

\begin{align}
\label{FFAforward}
F(k) = \sum_{l=-N/2}^{N/2-1} a_l e^{2\pi i \frac{l \cdot k}{2N-1}} , k:-(N-1),\cdots,N-1\\
\label{FFAbackward}
a_l = \frac{1}{2N-1} \sum_{k=-(N-1)}^{N-1} F(k) e^{-2\pi i \frac{l \cdot k}{2N-1}}, l : -\frac{N}{2},\cdots, \frac{N}{2}-1
\end{align}

Now if we look at the finite Fourier transform of $|F(k)|^2$ we have $|F(k)|^2 = \sum_{j=-(N-1)}^{N-1} b_j e^{2\pi i \frac{j \cdot k}{2N-1}}, k=-(N-1),\cdots,N-1$. On the other hand, $|F(k)|^2 = F(k)\cdot F^{\dagger}(k)$ with $F(k)$ expressed as in \eqref{FFAforward}. Those shall be equal, henceforth we have the following equation \eqref{FFArecursive}, written in the matrix format:


\begin{equation}
\label{FFArecursive}
 \underbrace{\begin{pmatrix}
  a^{\dagger}_{-\frac{N}{2}} & 0 & \cdots & 0 \\
  a^{\dagger}_{-\frac{N}{2} + 1} & a^{\dagger}_{-\frac{N}{2}} & \cdots & 0 \\
  \vdots  & \vdots  & \ddots & \vdots  \\
  a^{\dagger}_{\frac{N}{2}-1} & a^{\dagger}_{\frac{N}{2}-2} & \cdots & a^{\dagger}_{-\frac{N}{2}}
 \end{pmatrix}}_{\mathcal{A}}
\underbrace{
\begin{pmatrix}
a_{\frac{N}{2}-1}\\
a_{\frac{N}{2}-2}\\
\vdots \\
a_{-\frac{N}{2}}
\end{pmatrix}}_{\vec{a}}
=
\begin{pmatrix}
b_{N-1}\\
b_{N-2}\\
\vdots \\
b_{0}
\end{pmatrix}
\end{equation}

The phases of $a_l$ can be acquired from \eqref{FFArecursive} recursively by repeatedly solving equations of the type of \eqref{vector triangle} from the top line down, i.e. from the $b_{N-1}$ line to the $b_0$ line, assuming $N$ even. But at each step, we shall have a condition to choose a solution out of two. This is a discrete and nonlinear problem therefore a condition relying upon a `continuity' type argument might fail miserably. Here we adopt a closest point criteria. Namely, at the step of solving for  $a_{-N/2+l}$ and $a_{N/2-l-1}$, the two possible solutions are put in $\mathcal{A}$ and $\vec{a}$ as in \eqref{FFArecursive}, with all the calculated values of phases included. Unknown phases in $\mathcal{A}$ and $\vec{a}$ can be simply set as $1$, or some specific value if a priori is known. Then we compute the value of $\vec{c} = \mathcal{A} \vec{a}$ separately for the two possible solutions at this stage and choose the set of solution that makes the distance $\| \vec{c} - \vec{b} \|^2 = \|(c_N, c_{N-1}, \cdots, c_1) -  (b_N, b_{N-1}, \cdots, b_1)\|^2$ smaller. This is actually the same as the widely accepted idea of error reduction.

This closest point criteria can be proven to converge if the initial condition is self-consistent. The basic idea of the proof is like that of the contraction mapping principle. However, in reality the phase retrieval problem can be steep such that this criteria might not work exactly, i.e. there might be zeros in unwanted places, the function might have significant irregular high frequency components, and numerical accuracy issues, etc. Under those situations, this recursive method might only give an approximate solution given the nature of the closest point criteria.

\begin{figure}[htb]
\label{results1}
\centerline{\includegraphics[width=9cm]{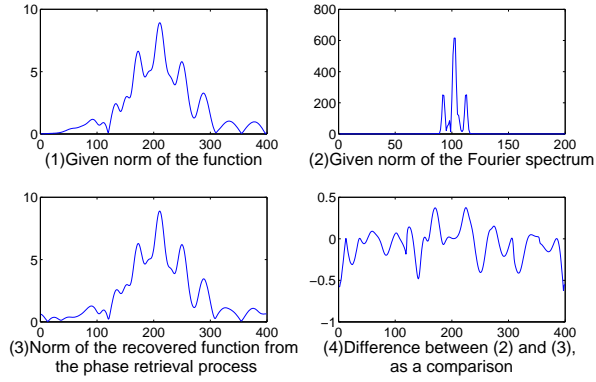}}
\caption{Result and behavior}
\end{figure}

The results of a algorithm based on this method is shown in Fig 2. The test function used is $h(t) = 2\exp( \frac{-(t-2)^2}{0.6} + 2it + 0.2(it-1)^2 )  (1.5 + 2\sin(5\pi t) + 0.5i\sin(2\pi t)+ \mathrm{sinc}(t-\pi)+3i\mathrm{sinc}(t-2)) + \sin(t^2)$ with 399 sampling points total for $t$ between $0$ and $3.98$ (with higher frequency components truncated). (1) and (2) in Fig 2 are the given conditions, i.e. norm of the function itself and its Fourier transform. Then the phases of the Fourier spectrum are recovered by the algorithm and an inverse Fourier transform is performed on the recovered spectrum to get back to the function itself. The norm of this computed result function is shown in (3) and a comparison is shown in (4).

The total number of calculations to solve this phase retrieval problem is $O(N^2)$, as we recursively solve $N/2$ triangles and at each of them it takes $O(N)$ calculations to choose the right one out of two.

\begin{figure}[htb]
\centerline{\includegraphics[width=9cm]{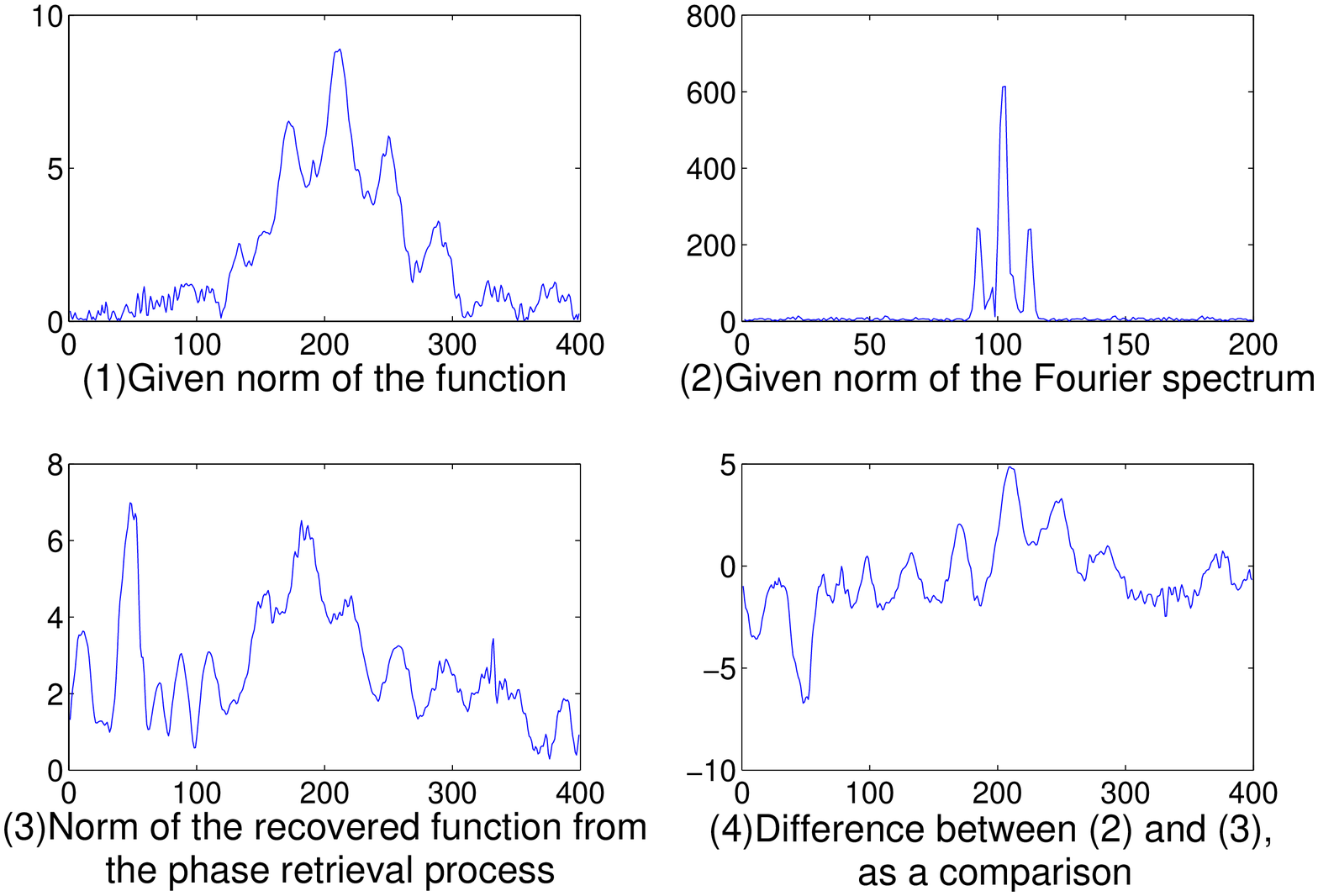}}
\caption{Behavior when subject to random noises}
\end{figure}

Since the recursive nature of this method, it is subject to the accumulation of numerical errors. Moreover, since it relies much upon the high frequency components, any instabilities or noises in the high frequency range could significantly reduce the performance. These are the major drawbacks of this elementary approach and it is hoped they are treated better in the future. A example with adding random noise is shown in Fig 3 where the test function is $h(t) + 0.3 \times \mathrm{random}(t)$.

\begin{figure}[htb]
\centerline{\includegraphics[width=9cm]{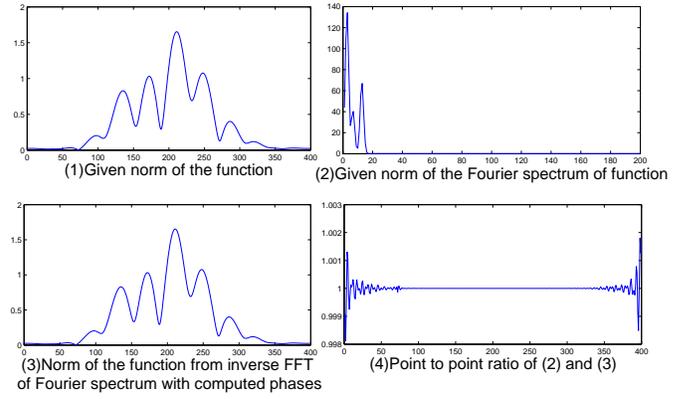}}
\caption{Behavior of this method when applied to functions without negative frequencies}
\end{figure}

This method can not be directly applied to the situation that both the norms of the function and its Fourier spectrum are given in terms of $N$ component vectors with no oversampling is possible. However, if more than half of the components in the vector of the Fourier spectrum are zero, this method has a natural extension to work under these situations. An extreme example of this is shown in Fig 4 where all the negative frequencies of the test function is zero. This elementary approach shall not be regarded as a final solution to the phase retrieval problem, and its performance on ill-defined phase retrieval problems still remains to be seen.

This recursive method can be combined with iterative methods. It would give an approximation instead of an exact solution, say, when it fails for a steep problem. That approximation can be used as a priori in the process of choosing a right solution when solving the triangles, with the algorithm operated from beginning once again. Or that approximation can serve as a starting point for an iterative method.

This method can as well be applied in higher dimensional phase retrieval problems. In a two dimensional problem under the context of Fourier series, suppose $f(x,y) = \sum_{n_1,n_2 = -N}^N a_{n_1 n_2}e^{in_1 x}e^{in_2 y}$. Then $ff^\dagger = \sum_{l_1,l_2 = -2N}^{2N} b_{l_1 l_2}e^{il_1 x}e^{il_2 y}$ with $b_{l_1 l_2} = \sum_{j1= max(l1-N,-N)}^{min(l1+N,N)} \sum_{j2= max(l2-N,-N)}^{min(l2+N,N)}a_{j_1 j_2} a^{\dagger}_{j1-l1 j2-l2}$. Then by starting with $b_{(2N)(2N)}$ and then move down step by step, the phases of $a_{n_1 n_2}$ could all recursively be recovered by solving equations of the type of \eqref{vector triangle}, combined with a similar closest point criteria. The difference now is the `counting geometry' in the recursive process. In the one dimensional problem the sequence follows $a_N, a_1 \rightarrow a_{N-1}, a_2 \rightarrow \cdots$. In the two dimensional problem, if the initial Fourier coefficients of $f$ are written in a matrix format, then this process begins with the coefficients on the very outside layer of the matrix, and then goes one layer by one layer further towards the core.



\end{document}